\documentclass[a4paper,11pt]{amsart}
\usepackage{amssymb,amsfonts,amsmath}
\usepackage{layout}
\usepackage[latin1]{inputenc}

\def\m{\mathcal}

\def\C{\mathbb{C}}
\def\c2{\mathbb{C}^2}
\def\R{\mathbb{R}}
\def\N{\mathbb{N}}

\def\N{\mathbb{N}}

\def\1{\bold{1}}

\def\a{\alpha}

\def\e{\varepsilon}

\def\f{\varphi}

\def\p{\psi}

\newcommand \W {\Omega}

\newcommand \mE {\mathcal E}
\newcommand \mF {\mathcal F}

\newcommand \vk {\chi}

\newcommand \sub {\subset}

\newcommand \mrm {\mathrm}

\newtheorem{lem}{Lemma}[section]
\newtheorem{pro}[lem]{Proposition}
\newtheorem{defi}[lem]{Definition}
\newtheorem{def/not}[lem]{Definition/Notations}
\numberwithin{equation}{section}
\newtheorem{thm}[lem]{Theorem}
\newtheorem{cor}[lem]{Corollary}

\newenvironment{proof3.1}
{\noindent {\it{Proof of theorem 3.1}}}{$\Box$ \linebreak[4]}

\begin{document}

\title[PSH functions  with weak singularities]
{Plurisubharmonic functions with weak singularities}

\author {S.BENELKOURCHI,  V.GUEDJ and A.ZERIAHI}
\maketitle
\begin{center}
{\it Dedicated to Professor C.O. Kiselman \\
on the occasion of his retirement}
\end{center}

\noindent \begin{abstract}
We study the complex Monge-Amp\`ere operator in bounded hyperconvex domains
of $\C^n$. We introduce several classes of weakly singular plurisubharmonic functions : these are functions  of finite
weighted Monge-Amp\`ere energy. They generalize the classes introduced by U.Cegrell,
and give a stratification of the space of (almost) all unbounded plurisubharmonic
functions. We give an interpretation of these classes in terms of the speed of decreasing
of the Monge-Amp\`ere capacity of sublevel sets and solve associated complex Monge-Amp\`ere equations.
\end{abstract}

{ 2000 Mathematics Subject Classification:} {\it 32W20, 32U05, 32U15}.

\section{Introduction}

In two seminal papers [Ce 1,2], U.Cegrell was able to define and study the
complex Monge-Amp\`ere operator $(dd^c \cdot)^n$ on special
classes of unbounded plurisubharmonic functions in a hyperconvex domain in $\C^n$.

Since we are considering a new and important scale of classes of plurisubharmonic functions  with finite weighted Monge-Amp\`ere energy, we find it convenient to introduce new notations which reflect our intuition. Therefore we have to modify some of the classical ones to avoid confusions.

Let $\W \subset \C^n$ be a bounded hyperconvex domain. 
The first important class considered by Cegrell (denoted by $\m E_0 (\W)$ in [Ce1]), is the class ${\mathcal T}(\Omega)$   of plurisubharmonic ``test functions'' on $\Omega$,
i.e. the convex cone of all bounded 
plurisubharmonic 
functions $\f$ defined on $\W$ such that $\lim _{z\to \zeta} \f (z) = 0,$ 
for every $\zeta \in \partial \W,$ and $\int_\W(dd^c \f )^n <+\infty$.
Besides this class, we will need the following classes introduced in [Ce1], [Ce2].
\begin{itemize}
\item The class $DMA(\Omega)$ 
 is the set of plurisubharmonic functions
 $u $
such that for all $z_0 \in \Omega $,
 there exists a neighborhood $V_{z_0}$ of
$z_0$ and $u_j \in {\mathcal T}(\Omega)$ a decreasing sequence which
converges towards $u$ in $V_{z_0}$ and satisfies
$\sup_j \int_{\Omega} (dd^c u_j)^n <+\infty$.
U.Cegrell has shown [Ce 2] that the operator $(dd^c \cdot )^n$ is well defined
on $DMA(\Omega)$ and continuous under decreasing limits. The class
$DMA(\Omega)$ is stable under taking maximum and it is the largest
class with these properties (Theorem 4.5 in [Ce 2]).
Actually this class, introduced and denoted by $\m E (\Omega)$  by U.Cegrell ([Ce 2]),  turns out  to  coincide with the domain of definition of the complex Monge-Amp\`ere operator on $\Omega$ as was shown by Z.Blocki [Bl 1,2];

\item the class ${\mathcal F}(\Omega)$ is the ``global version'' of $DMA(\Omega)$:
a function $u$ belongs to ${\mathcal F}(\Omega)$ iff there exists $u_j \in {\mathcal T}(\Omega)$
a sequence decreasing towards $u$ {\it in all of } $\Omega$, which satisfies
$\sup_j \int_{\Omega} (dd^c u_j)^n<+\infty$;
\item the class ${\mathcal F}_a (\Omega)$ is the set of functions $u \in {\mathcal F}(\Omega)$
whose Monge-Amp\`ere measure $(dd^c u)^n$ is absolutely continuous with respect to capacity i.e. it does not charge pluripolar sets;
\item the class ${\mathcal E}^p(\Omega)$ (respectively  ${\mathcal F}^p(\Omega)$) is the
set of functions $u$ for which there exists a sequence of
functions $u_j \in {\mathcal T}(\Omega)$ decreasing towards $u$ in all of $\Omega$, and
so that $\sup_j \int_{\Omega} (-u_j)^p (dd^c u_j)^n<+\infty$
(respectively $\sup_j \int_{\Omega} [1+(-u_j)^p] (dd^c u_j)^n<+\infty$).
\end{itemize}
One purpose of this article is to 
use the formalism developed in [GZ] in a compact setting
to give a unified treatment of all these classes.
Given an increasing function $\chi:\R^- \rightarrow \R^-$, we consider the
set ${\mathcal E}_{\chi}(\Omega)$ of plurisubharmonic functions
of finite $\chi$-weighted Monge-Amp\`ere energy. These are functions
$u \in PSH(\Omega)$ such that there exists $u_j \in \m T(\W)$ decreasing to $u$, with
$$
\sup_{j \in \N} \int_{\Omega} (-\chi) \circ u_j (dd^c u_j)^n <+\infty.
$$
It will be shown that ${\mathcal E}_{\chi}(\Omega) \sub DMA (\Omega)$.

 Many important properties
follow from the elementary observation  that the Monge-Amp\`ere 
measures $ \mrm 1 _{\{u>-j\}}  (dd^c u_j)^n$ strongly converge towards $(dd^c u)^n$
in the set  $\Omega \setminus (u=-\infty)$, when  $u_j := \max (u , - j)$ are the "canonical approximants'' of $u$:
\vskip.2cm
\noindent {\bf Theorem A.} {\it
If $u\in DMA(\W ),$ then  for all Borel sets 
 $B\subset \W \setminus \{u=-\infty\}$,
$$
\int_B(dd^c u )^n 
 = \lim _{j\to \infty } \int_{ B\cap \{u>-j\}}  (dd^c u_j )^n  ,
$$
where $u_j :  = \max (u, -j )$ are the canonical approximants.
}
\vskip.2cm

We establish this result in {\it section 2} and derive several consequences.
This yields in particular simple proofs of quite general comparison principles.

\vskip.1 cm
The classes ${\mathcal E}_{\chi}(\Omega)$ have very different properties, depending
on whether $\chi(0)=0$ or $\chi(0) \neq 0$, $\chi(-\infty)=-\infty$ or 
$\chi(-\infty) \neq -\infty$, $\chi$ is convex or concave.
We study these in {\it section 3} and give a capacitary interpretation of them
in {\it section 4}. Let us stress in particular Corollary 4.3 which gives
an interesting characterization of the class ${\mathcal E}^p(\Omega)$
of U.Cegrell, in terms of the speed of decreasing
of the capacity of sublevel sets:
\vskip.1cm

\noindent {\bf Proposition B.} {\it For any real number $p > 0$,
$$
\mathcal E^p (\Omega) = \left\{ \f \in PSH^- (\Omega) ; \int_0^{+ \infty} (- \f)^{n + p - 1} 
Cap_{\Omega} (\{\f < - t\}) d t < + \infty \right\}.
$$
}

Here $Cap_{\Omega}$ denotes the Monge-Amp\`ere capacity introduced by E. Bedford and B.A. Taylor ([BT1]). Of course $\mathcal E^p (\Omega) = \mathcal E_{\chi} (\Omega),$ for $\chi (t) := - (- t)^p$.
\vskip.2cm

Our formalism allows us to consider
further natural subclasses of $PSH(\Omega)$, especially functions
with finite ``high-energy'' (when $\chi$ increases faster than polynomials
at infinity). We study in {\it section 5} the range of the Monge-Amp\`ere 
operator on these classes.
Given a positive finite Borel measure $\mu $ on $\Omega$, we set
$$F_{\mu} (t) := \sup \{ \mu (K) ; K \subset \W  \  \ \mbox{compact,  } \ \ Cap_{\Omega} (K) \leq t\}, t \geq 0.$$
 Observe that $ F := F_{\mu}$ is an increasing function on $\mathbb R^{+}$ which satisfies
$$
\mu(K ) \le F (Cap_{\Omega} (K) ), \quad \mbox{for all Borel subsets } \ K \subset X.
$$
The measure $\mu$ does not charge pluripolar sets iff $F (0) = 0$.

When $F (x) \lesssim x^\alpha $ vanishes 
at order 
$\alpha >1$, S. Kolodziej has proved \cite{K 2} that the equation
$\mu=(dd^c \f)^n$ admits a unique {\it continuous} solution 
with $\f_{|\partial \Omega}=0$.
If $F(x) \lesssim x^\alpha $  with $0< \alpha <1 ,$  
it follows from the work of
U. Cegrell \cite{Ce 1} that 
there is a unique solution in some class ${\mathcal F}^p(\Omega)$.

Another objective of this article is to fill in the gap inbetween Cegrell's and 
Kolodziej's results, by considering all intermediate dominating functions $F.$ 
Write $ F(x) = x [\e( -\ln x /n)]^n$ where  $\e:\R^+ \rightarrow  [0 ,\infty [$ 
is nonincreasing. 

Our second main result is:

\vskip.2cm
\noindent {\bf Theorem C.} {\it 
Assume for all compact subsets $K  \subset \W$,
$$
\mu(K) \leq F_{\e} (\mrm{Cap}_{\Omega}(K)) , 
\text{ where } F_{\e}(x)=x[\e(-\ln x/n)]^n.
$$
Then there exists a unique function $\f \in \mF (\W)$ 
such that $\mu=(dd^c \f)^n$
and
\begin{equation*}
\mrm{Cap}_{\Omega}(\{\f<-s \}) \leq  \exp (-nH^{-1}(s)), \text{ for all }\ s>0,
\end{equation*}
Here $H^{-1}$ is the reciprocal function of
 $H(x) =e\int_{0}^x  \e(t) dt +s_0(\mu)$.

In particular $\f \in \mE _\chi(\Omega) $ where
$-\chi(-t)=\exp (  n H^{-1}(t)/2)$.}
\vskip.2cm
 Note in particular that when $\mu \leq Cap_{\Omega}$ (i.e. $\e \equiv 1$), then $\mu = (dd^c \f)^n$ for a function $\f \in \mF (\Omega)$ such that $Cap_{\Omega} (\{\f < - s\})$ decreases exponentially fast. Simple examples show that this bound is sharp (see [BGZ]).

For similar results in the case of compact K\"ahler manifolds, we refer the reader
to [GZ], [EGZ], [BGZ].

\vskip.2cm
\noindent{\bf Remerciements.} {\it C'est un plaisir de contribuer \`a ce volume en l'honneur de Christer Kiselman, dont nous avons toujours appr\'eci\'e la gentillesse et la grande  \'el\'egance math\'ematique.}

\section{Canonical approximants}

We let $PSH(\W) $ denote the set of plurisubharmonic functions on $\W$ (psh for short),
and fix $u \in PSH(\Omega)$.
E.Bedford and B.A.Taylor have  defined in \cite{BT 2}  the non pluripolar part
 of the Monge-Ampère
measure of $u$: the sequence    
 $ \mu^{(j)}_u: = \mrm{1}_{\{u>-j\}} (dd^c \max[u,-j] )^n $ is a  nondecreasing sequence of positive measures.
Its limit $\mu_u $ is the ``nonpluripolar part of $(dd^c u )^n$'', defined as,  
$$
\mu_u (B) = \lim _{j \to \infty }\int _{B\cap \{u>-j\}} (dd^c \max [u, -j])^n,
$$
for any Borel set $B\subset \W.$

 In general $\mu_u $ is not  locally bounded near $\{u=-\infty \}$ (see e.g. \cite{Ki}),
 but if $u\in DMA(\W) $ then
 $\mu_u $ is a regular Borel measure:

\begin{thm} \label{nonplp}
If $u\in DMA(\W ),$ then  for all Borel sets 
 $B\subset \W \setminus \{u=-\infty\}$,
$$
\int_B(dd^c u )^n 
 = \lim _{j\to \infty } \int_{ B\cap \{u>-j\}}  (dd^c u_j )^n  ,
$$
where $u_j :  = \max (u, -j ).$
In particular, $\mu _u =\mrm{1}_{\{u>-\infty \}}(dd^c u )^n$.

The measure $(dd^c u )^n $ puts no mass on pluripolar sets $E \subset \{u >-\infty \}.$
\end{thm}

\begin{proof} 
Note that this convergence result is local in nature,
 hence we can assume, without loss of generality,
 that $u \in \mF  (\W).$
For $s>0$ consider the psh  function $h_s : = \max (u/s +1 ,  0)$. Observe that $h_s$
 increases to the Borel function $\mrm{1}_{\{u>-\infty\}}$ and
 $ \{h_s =0\} = \{u\le -s\}$. We claim that 
$$
h_s (dd^c  \max (u, -s ))^n = h_s (dd^c u)^n \ , \  \mbox{for all }\ s>0,
$$
in the sense of measures on $\W$.

Indeed, recall that we can find a sequence of {\it continuous} tests functions
 $u_k $ in  $\m T(\W) $ decreasing towards $u$  (see Theorem 2.1 in \cite{Ce 2}). 
It follows from Proposition 5.1 in  \cite{Ce 2} that 
$h_s (dd^c  \max (u_k , -s ))^n $ converges weakly to $ h_s (dd^c  \max (u, -s ))^n$  and 
$ h_s (dd^c u_k)^n $ converges weakly to $ h_s (dd^c u)^n$ as $k\to \infty $.

 Since  
$ \max (u_k , -s ) = u_k $  on $ \{u_k  >-s\}$, which is an open neighborhood of the set
 $\{u>-s\},$
we infer  $$ h_s (dd^c  \max (u , -s ))^n =  h_s (dd^c u)^n, $$
 as claimed. 

Observe that
 $$
h_s (dd^c  \max (u, -s ))^n = h_s \mrm{1}_{\{u>-s\}} (dd^c u)^n = h_s \mu_u^{(s)}
 $$
increases as $s \uparrow + \infty$
towards $\mrm{1}_{\{u>-\infty \}} \mu_u = \mu_u $, as follows from the monotone convergence  and Radon-Nikodym theorems. Similarly 
$ h_s (dd^c u)^n$  converges to  
$ \mrm{1}_{\{u>-\infty\}} (dd^c  u )^n. $  Thus  
  $\mu_u = \mrm{1}_{\{u> -\infty\}} (dd^c  u )^n, $ this shows the desired convergence on any Borel set $B\subset \W \setminus \{u= - \infty\}$.
\end{proof}

Note that if $u\in \mF_a (\W) $ then 
$
\int_B(dd^c u )^n 
 = \lim _{j\to \infty } \int_{ B}  (dd^c u_j )^n  ,
$
for all Borel subsets $B\subset \W$ (see Theorem 3.4).

As an application, we give a simple proof of the following general
version of the comparison principle (see also \cite{NP}).

\begin{thm}\label{demailly}
Let $u  \in DMA (\W)$ and $v\in PSH^-(\W).$ Then 
$$
{\bf 1}_{\{u>v\}} (dd^c u )^n ={\bf 1}_{\{u>v\}} (dd^c \max ( u , v) )^n 
$$
\end{thm}

\begin{proof}
 Set $u_j = \max (u, -j ) $ 
and $v_j =\max (v, -j )$. Recall from \cite{BT 2} that the desired equality is
 known for bounded psh functions, 
$$
{\bf 1}_{\{u_j>v_{j+1}\}} (dd^c u_j )^n ={\bf 1}_{\{u_j>v_{j+1}\}}
 (dd^c \max ( u_j , v _{j+1}) )^n .
$$
Observe that 
 $\{u>v \} \subset \{u_j>v_{j+1}\}$, hence 
\begin{multline*}
{\bf 1}_{\{u>v \}}\cdot {\bf 1}_ {\{u  >-j \} } (dd^c u_j )^n =
{\bf 1}_{\{u>v \}}\cdot {\bf 1}_ {\{u  >-j \} }  (dd^c \max ( u , v , -j) )^n 
\\
 =
{\bf 1}_{\{u>v\}} \cdot {\bf 1}_{ \{\max (u , v) >-j \} } (dd^c \max ( u , v , -j) )^n 
.
\end{multline*}
It follows from  Theorem \ref{nonplp} that
$ {\bf 1}_ {\{u  >-j \} } (dd^c u_j )^n $ converges in the strong sense of Borel 
measures towards $\mu_u = {\bf 1}_{\{ u>-\infty\}}(dd^c u )^n  $.
 Observe that  ${\bf 1}_{\{u>v\}} {\bf 1}_{\{ u>-\infty\}} = {\bf 1}_{\{u>v\}} $. 
We infer, by using Theorem \ref{nonplp} again with $\max ( u , v)$, that
$$
{\bf 1}_{\{u>v\}} (dd^c u )^n ={\bf 1}_{\{u>v\}} (dd^c \max ( u , v) )^n .
$$ \end{proof}
The following result has been  proved by U.Cegrell \cite{Ce 3}. We provide here a 
simple proof using Theorem 2.2, yet another consequence of the fact that the Monge-Amp\`ere measures $\mrm 1 _{\{u >-j\}} (dd^c u_j)^n$ strongly converge towards
 $\mrm 1 _{\{u >-\infty\}}(dd^c u)^n$ when $u_j:=\max(u,-j)$
are the ``canonical approximants'' (Theorem \ref{nonplp}).

\begin{cor}\label{cegrell3}
Let $\f \in \mF(\W)$ and $u \in DMA (\W)$ such that $u \leq 0$.  Then 
$$ 
\int _{\{  \f < u \} } (dd^c u )^n \le \int _{\{  \f < u \} \cup \{\f = -\infty\} } (dd^c \f )^n 
$$
\end{cor}
\begin{proof}
Since $\psi := \max \{u, \f \} \in  \mF(\W)$ and $\f \leq \psi$ on $\W$, it follows  that 
 $$
 \int _\W (dd^c \psi)^n   \le 
 \int _\W (dd^c \f )^n.
$$
Indeed this is  clear when $\f \in \mathcal T (\W)$ by integration by parts and follows by approximation when $\f \in \mathcal F (\W)$ (see \cite{Ce 2}).

We infer by using Theorem \ref{demailly}, 
\begin{eqnarray*}
\int _{\{\f < u\}} (dd^c u )^n &=& \int _{\{\f < u\}} (dd^c \max(u,\f) )^n \\
&=& \int _{\W} (dd^c \max(u,\f) )^n - \int _{\{\f \ge u\}} (dd^c \max(u,\f) )^n \\
&\le &  \int _{\W} (dd^c \f )^n - \int _{\{\f > u \}} (dd^c \f )^n  -
\int _{\{\f = u \}} (dd^c \max(u,\f) )^n\\
&\le & \int _{\{\f \le u\}} (dd^c \f)^n
\end{eqnarray*}
Now take  $0< \e <1$ and apply the previous result to get
$$
\int _{\{ \e  \f < u \} } (dd^c u )^n \le \int _{\{  \e \f \le u \} } (dd^c \e \f  )^n
=
\e ^n  \int _{\{  \e \f \le u \} } (dd^c \f   )^n.
$$
The desired inequality follows by  letting $\e \to 1,$ since  $ \{  \e \f < u \} $ increases to
 $\{  \f < u \}$
and 
$\{  \e \f \le u\} $  increases to $\{  \f < u \}  \cup \{\f = -\infty\}.$

\end{proof}

Note that Corollary \ref{cegrell3} is still valid when $\f , u \in DMA(\W)$ under the condition 
 $\{\f < u\} \Subset \W.$

The following comparison principle is due to  U.Cegrell (see 
Theorem 5.15 in \cite{Ce 2} and Theorem 3.7 in \cite{Ce 3}). 

\begin{cor} \label{cegrell}
Let  $\f   \in \mF_a (\W ) $ and $ u \in DMA(\W),$  
such  that 
$(dd^c \f )^n \le  (dd^c u )^n .$  Then  $u \le \f .$

In particular if $(dd^c u )^n =  (dd^c \f )^n $ with 
$u , \f  \in \mF_a (\W )$, then $u=\f.$
\end{cor}
\begin{proof} The proof is a consequence of Corollary \ref{cegrell3} and follows from standard arguments  (see e.g.  \cite{BT 1} for bounded psh function).
\end{proof}
Note that the result still holds when   $u \in DMA (\W)$ is  such that $(dd^c u)^n $ vanishes on pluripolar sets and $u\ge v $ near
 $\partial \W$. 
However it   fails in $\m F (\W)$ (see \cite{Ce 2} and \cite{Z}).

Now, as another consequence of  Theorem 2.2, we provide the following result which will be
 useful in the sequel:
\begin{cor}\label{est}
Fix $\f \in {\mathcal F}(\W)$. Then  for all $s>0$ and $ t > 0$,
\begin{equation} \label {eq:FI}
t^n Cap_\W(\{\f < -s-t\} ) \leq \int_{(\f <  -s)} (dd^c \f)^n \leq s^n Cap_\W(\{\f <  -s\} ).
\end{equation}
In particular 
\begin{equation} \label{FE}
 \int_\W(dd^c \f )^n = \lim _{s\downarrow 0} s^n Cap_\W (\ \le -s) = \sup _{s > 0} s^n Cap_\W (\f < -s).
\end{equation}
Moreover a negative function $ u \in PSH (\W)$ belongs to ${\mathcal F}(\W)$ if and only if 
$ \sup _{s > 0} s^n Cap_\W (u < -s) < + \infty$
\end{cor}
The inequalities  (\ref{eq:FI}) was proved for psh test functions in [K3] (see also [CKZ] and [EGZ]).  For $\f \in {\mathcal F}(\W) \cap L^{\infty} (\Omega),$ it follows by approximation and quasi-continuity. In the general case, it can be deduced using Theorem 2.1.  The last assertion follows easily from (\ref{eq:FI}). It was first obtained in ([B]).

\section{Weighted energy classes}

\begin{defi}
 Let $\chi : \R^-  \to \R^- $ be an
 increasing
 function. We let $ \m E_\chi (\W ) $
 denote the set of all functions $u \in PSH(\Omega)$ for which
 there exists a sequence $u_j \in {\mathcal T}(\Omega)$
 decreasing to $u$ in $\Omega$ and satisfying
 $$ 
\sup_{j \in \N} \int _ {\W}(-\chi) \circ u_j \, (dd^c u_j )^n <\infty.
$$
\end{defi}

This definition clearly contains the  classes of U.Cegrell:
\begin{itemize}
\item $\m E_\chi (\W )={\mathcal F}(\W)$ if $\chi$ is bounded and $\chi(0) \neq 0$;
\item $\m E_\chi (\W )={\mathcal E}^p(\W)$ if $\chi(t)=-(-t)^p$;
\item $\m E_\chi (\W)={\mathcal F}^p(\W)$ if $\chi(t)=-1-(-t)^p$.
\end{itemize}

We will give hereafter interpretation of the classes 
   ${\mathcal F}(\Omega) \cap L^{\infty}(\Omega)$ and ${\mathcal F}_a(\Omega)$
in terms of weighted-energy as well.

Let us stress that the classes 
${\mathcal E}_{\chi}(\Omega)$ are very different whether $\chi (0) \neq 0 $ 
(finite total Monge-Amp\`ere mass) or  $\chi (0)=0$.

To simplify we consider in this section the case $\chi(0) \neq 0$, so that
all functions under consideration have a well defined Monge-Amp\`ere
measure of finite total mass in $\Omega$.  Note however that many results
 to follow still hold when $\chi(0)=0.$

\begin{pro}\label{prop}
Let $\chi: \R^- \rightarrow \R^-$ be an increasing
function such that $\chi(-\infty)=-\infty$ and $\chi(0) \not = 0$. Then
$$
{\mathcal E}_{\chi}(\Omega) \subset {\mathcal F}_{a}(\Omega).
$$
In particular the Monge-Amp\`ere measure $(dd^c u)^n$ of a function
 $u \in \mE _\chi (\W )$
is well defined and does not charge pluripolar sets. More precisely,
$$
\mE_\chi (\W)=\left\{ u \in {\mathcal F}(\Omega) \, / \, 
\chi \circ u \in L^1((dd^c u)^n) \right\}.
$$\end{pro}

\begin{proof} Fix $u\in \mE_\chi (\W) $ and $u_j \in \m T(\W )$ a defining sequence 
 such that 
$$\sup _j \int_\W \chi (u_j)(dd^c u_j )^n <+\infty .$$
 The condition $\chi(0) \not = 0$ implies that 
$
{\mathcal E}_{\chi}(\Omega) \subset {\mathcal F}(\Omega).
$
In particular the Monge-Ampère measure $(dd^c u )^n $ is well defined. It follows from the upper semi-continuity of $u$ that $ - \chi (u)(dd^c u )^n $ is bounded from above by any cluster point of the bounded sequence $- \chi (u_j)(dd^c u_j )^n.$ Therefore 
$\int_\W (-\chi )\circ u (dd^c u )^n < +\infty , 
$  in particular $(dd^c u )^n $ does not charge the set $\{\chi(u) = -\infty \},$ which
 coincides with $\{u= -\infty \}$, since $\chi(-\infty ) = -\infty .$ It follows therefore from   Theorem 2.1 that the measure
$ (dd^c u)^n$ does not charge pluripolar sets.

To prove the last assertion, it remains to show the reverse inclusion  $$
\mE_\chi (\W)\supset\left\{ u \in {\mathcal F}(\Omega) \, / \, 
\chi \circ u \in L^1((dd^c u)^n) \right\}.$$ 
So fix $u\in \mF(\W )$  such that $\chi \circ u \in L^1((dd^c u)^n).$
It follows from \cite{K 1} that  there exists, for each $ j\in \N $,  a function 
$u_j \in \mathcal{T} (\W ) $ 
 such that $(dd^c u_j )^n = {\bf 1}_{\{u>j\rho \}} (dd^c u)^n$, where $\rho \in \m T (\W)$
any defining function for $\W = \{\rho < 0\}$.
Observe  that $ (dd^c u )^n \ge (dd^c u_{j+1} )^n \ge  (dd^c u_j )^n$. We infer
  from Corollary \ref{cegrell} that  $(u_j) $  is a decreasing sequence 
and   $u\le u_j$. The monotone convergence theorem thus yields 
$$
\int_\W  (- \chi)\circ  u_j (dd^c u_j )^n = \int _\W (- \chi) \circ u_j {\bf 1}_{\{u>j\rho \}} 
(dd^c u)^n \to  
\int_\W (-  \chi ) \circ  u (dd^c u )^n <+\infty,
$$
so that $u\in \mE_\chi (\W).$
\end{proof}

There is a natural partial ordering of the classes $\mE_{\chi} (\W):$
if $ \chi = O (\tilde{\chi })$ then $\m E_ {\tilde{\chi }}(\W) \subset \m E_ {{\chi }} (\W)$. 
 Classes ${\mathcal E}_{\chi}(\Omega)$ provide a 
full scale of subclasses of $PSH^-(\W )$  
of unbounded functions, reaching, 
``at the limit'', bounded plurisubharmonic functions.

\begin{pro}
$$
\quad   \quad \mF(\W )\cap L^\infty(\Omega)  =
 \bigcap_{\substack{ \vk(0)\not =0 \\
\chi(-\infty )=-\infty} } \mE_\vk (\W ),
 %\quad \mF = \bigcap_{ \vk\  bounded} \mE_\vk .
$$
where the intersection runs over all increasing functions 
$\chi:\R^- \rightarrow \R^-$.
\end{pro}

Note that it suffices to consider here those functions $\chi$ which are {\it concave}.

\begin{proof}
One inclusion is clear. Namely  if $u\in \mF(\W )\cap L^\infty(\Omega)$ and $u_j \in \m T (\W) $ are decreasing to $u,$ then for any  $\chi $ as above,
$$
\int _\W -\chi (u_j) (dd^c u_j )^n \le \left [\sup _\W |\chi (u)|\right ] \int _\W (dd^c u)^n <+\infty.
$$
Conversely, assume  $u\in 
 \mF (\W ) $ is unbounded. Then the sublevel sets $\{u <t \} $ 
are non empty for all $t<0,$ hence 
 we can  consider the function   $\chi $ such that  
 $$
t \mapsto \chi^\prime (t) = \frac{1}{(dd^c u)^n(\{u<t\})}, \mbox {for all } t<0.
$$
The function  $\chi $ is clearly  increasing. Moreover $ (dd^cu )^n$  has finite (positive) mass,
 hence $\chi^{\prime} (t) \ge \frac{1}{(dd^cu )^n (\W) }.$ This yields  
$\chi(-\infty )=-\infty$.  Now   
$$
\int _\W (-\chi) \circ u  (dd^c u )^n = \int _0 ^{+\infty}\chi^\prime (-s)(dd^c u)^n(\{u<-s\})ds = + \infty. 
$$
 This  shows that if $u\in \m E_\chi (\W )$ for all $\chi $ as above, then $u$ has to be bounded.
\end{proof}

When $u\in \mE_\chi (\W ) \subset \mF_a(\W) $, the canonical approximants
  $u_j : =  \max (u, -j ) $ yield strong convergence properties
 of weighted  Monge-Ampère operators:

\begin{thm}
Let $\chi: \R^- \rightarrow \R^-$ be an increasing
function such that $\chi(-\infty)=-\infty$ and $\chi(0) \not = 0$. Fix $u\in \mE_\chi(\W)$ as set  $u^j =  \max (u , -j).$
Then for each Borel subset $B\subset \W ,$
$$
\lim_{j \to + \infty}  \int_B \chi (u^j)(dd^c u^j)^n = \int_B \chi (u)(dd^c u)^n. $$
Moreover if   $(u_j)_{j \in \N}$ is any decreasing sequence in $\mathcal E_{\chi} (\Omega)$  converging to $u$ such that $\sup_j \int_{\W} \vert \chi (u_j)\vert  (dd^c u_j)^n < + \infty$, then 
$$ 
\lim_{j \to + \infty}  \int_{\W} \chi (u_j) (dd^c u_j)^n = \int_{\W} \chi (u) (dd^c u)^n. 
$$
\end{thm}

Let us stress that this convergence result is stronger than Theorem 5.6 in [Ce 1]:
on one hand we produce here an explicit (and canonical) sequence
of bounded approximants, on the other hand the convergence holds
in the strong sense of Borel measures. Moreover the $\chi-$energy is continuous under  decreasing sequences of plurisubharmonic functions with uniformly bounded $\chi-$energies.

\begin{proof} 
We first show that $(dd^c u^j )^n $ converges towards $(dd^c u )^n $ ``in the strong sense 
of Borel measures'', i.e. $(dd^c u^j )^n (B) \to (dd^c u )^n (B), $ 
for any Borel set $B\subset \W .$  
Observe that for $j\in \N^* $ fixed and $0 < s < j $, 
$\{ u<-s \} = \{ u_j < -s\}$. 
It follows from Corollary  \ref{est} that 
%\cite{B} (Proposition 2.1 ) that 
 $$\int_\W (dd^c u^j)^n= \int_\W (dd^c u)^n .$$
Therefore 
\begin{eqnarray*}
\int_{\{u\le -j\}} (dd^c u^j)^n & =&\int_\W (dd^c u^j)^n  - \int_{\{u> -j\}} (dd^c u^j)^n\\ 
&=& \int_\W (dd^c u)^n  - \int_{\{u> -j\}} (dd^c u)^n =  \int_{\{u\le -j\}} (dd^c u)^n . 
\end{eqnarray*}
 Thus if $B \subset \W $ is a Borel subset,
\begin{eqnarray*}
\left |\int _B (dd^c u^j) ^n -\int _B (dd^c u )^n \right | &\le & \int _{ \{u\le -j\}}
 (dd^c u^j) ^n + 
\int _{\{u\le -j\}} (dd^c u )^n \\
&\le & 2 \int _{\{u\le -j\}} (dd^c u )^n \to 0 ,\  \mbox{as} \ j\to + \infty.
\end{eqnarray*}
The proof that $\chi \circ u^j (dd^c u^j )^n $ converges strongly towards
 $\chi \circ u (dd^c u )^n $  goes along  similar lines, once we observe that
 \begin{multline*}
\int_{\{u\le -j \}} -\chi \circ u^j (dd^c u^j )^n  = -\chi (-j )\int_{\{u\le -j \}} 
 (dd^c u^j )^n = \\ -\chi (-j )\int_{\{u\le -j \}} 
 (dd^c u )^n \le \int_{\{u\le -j \}} -\chi \circ u  (dd^c u )^n .
\end{multline*}

To prove the second statment we proceed as in [GZ].  Observe that the statement is true for uniformly bounded sequences of plurisubhatmonic functions by Bedford and Taylor convergence theorems.  For the general case, we first consider  an increasing function $\tilde \chi : \R^- \longrightarrow \R^-$ such that $\tilde \chi = o (\chi)$ and prove the convergence of the $\tilde \chi-$energies.
Indeed, for $k \in \N$ define the canonical approximants
$$
u_j^{k} := \sup \{u_j,-k\}, \ \ \mathrm{and} \quad u^{k} := \sup \{u,-k\}.
$$
The integer $k$ being fixed, the sequence $(u_j^k)_{j \in \N}$ is uniformly bounded and decreases towards $u^k$, hence the $\tilde \chi-$energies of $ u_j^k$ converge to the $\tilde \chi-$energy of $u^k$ as $j \to + \infty$.
Thus we will be done if we can show that  the $\tilde \chi-$energies of $ u_j^k$ converge to the $\tilde \chi-$energy of $u_j$ uniformly in  $j$ as $k \to + \infty$. This follows easily from the following inequalities
\begin{eqnarray*} 
I (j,k) &:=&  \left \vert \int_{\W} \tilde \chi (u_j^k)(dd^c u_j^k)^n -  \int_{\W} \tilde \chi (u_j)(dd^c u_j)^n \right\vert \\
&\leq &  \int_{\{u_j \leq - k\}} - \tilde \chi (u_j^k)(dd^c u_j^k)^n  + \int_{\{u_j \leq - k\}} - \tilde \chi (u_j)(dd^c u_j^k)^n \\
&\leq &  \frac{\tilde \chi (-k)}{ \chi (-k)} \left(\int_{\{u_j \leq - k\}} -  \chi (u_j^k)(dd^c u_j)^n  + \int_{\{u_j \leq - k\}} -  \chi (u_j)(dd^c u_j)^n\right) \\
& \leq & 2 \frac{\tilde \chi (-k)}{ \chi (-k)} \int_{\W} - \chi (u_j)(dd^c u_j)^n \leq  2 M \frac{\tilde \chi (-k)}{ \chi (-k)},
\end{eqnarray*}
where $M := \sup_j \int_{\W} - \chi (u_j)(dd^c u_j)^n < + \infty$ and the last inequality follows from previous computations.

For the general case, observe that  $ 0 \leq f := - \chi (u) \in L^1 ((dd^c u)^n)$ by Proposition 3.2. Then  it follows easily by an elementary integration theory argument that there exists an increasing function $h : \R^+ \longrightarrow \R^+$ such that $\lim_{t \to + \infty} h (t) \slash t = + \infty$ and $ h (f) \in  L^1 ((dd^c u)^n)$ (see [RR]). Thus  $u \in \mathcal E_{\chi_1} (\W)$, where $\chi_1 (t) := - h (- \chi (t))$ for $t < 0$ and  $\chi = o (\chi_1)$ and  the continuity property for $\chi-$energies follows from the previous case.

\end{proof}

\section{Capacity estimates}

Of particular interest for us here are the classes
${\mathcal E}_{\chi}(\Omega)$, where the weight $\chi:\R^- \rightarrow \R^-$
has fast growth at infinity.
It is useful in practice to understand these classes through the speed of decreasing of the capacity
of sublevel sets.

The Monge-Amp\`ere capacity has been introduced and studied by E.Bedford and
A.Taylor in [BT 1]. Given $K \subset \W$ a Borel subset,
it is defined as
$$
\mrm {Cap }_\W(K):=\sup \left\{ \int_K (dd^c u)^n \, / \, u
\in PSH(\W), -1\leq u \leq 0 \right\}.
$$

\begin{defi}
$$
\hat{{\mathcal E}}_{\chi}(\W ) :=\left\{ \f \in PSH(\W) \, / \,
\int_0^{+\infty} t^n \chi'(-t) \mrm {Cap}_\W(\{\f<-t\})
dt<+\infty \right\}.
$$
\end{defi}

The classes ${\mathcal E}_{\chi}(\Omega)$ and $\hat{{\mathcal E}}_{\chi}(\Omega)$
are closely related:

\begin{pro}
The classes $\hat{{\mathcal E}}_\chi(\W ) $ are convex and stable under maximum:  if $\f \in \hat{{\mathcal E}}_{\chi}(\W ) $ and $  \p \in PSH^-(\W)$,
then $\max(\f, \p )  \in \hat{{\mathcal E}}_{\chi}(\W)$. 

One always has
$\hat{{\mathcal E}}_{\chi} (\W ) \subset \mE _\chi (\W ) $, while
$$
{\mathcal E}_{\hat{\chi}}(\W) \subset \hat{{\mathcal E}}_{\chi}(\W), \text{
where }  \hat{\chi}(t) = \chi(2t).
$$
\end{pro}
 
\begin{proof}
The convexity of $\hat{{\mathcal E}}_{\chi}(\W)$ follows from
the following simple observation: if $\f,\p \in \hat{{\mathcal E}}_{\chi}(\W)$
and $0 \leq a \leq 1$, then
$$
\left\{ a\f+(1-a)\p <-t \right\} \subset
\left\{ \f<-t \right\} \cup  \left\{ \p <-t \right\}.
$$
The stability under maximum is obvious.

Assume $\f \in \hat{{\mathcal E}}_{\chi}(\W)$.
We can assume without loss of generality $\f \leq 0$ and $\chi(0)=0$.
Set $\f_j:=\max(\f,-j)$. It follows from Corollary \ref{est} that
\begin{eqnarray*}
\int_\W (-\chi) \circ \f_j \, (dd^c {\f_j})^n &=&
\int_0^{+\infty} \chi'(-t) (dd^c {\f_j})^n(\f_j < -t) dt \\
&\leq& \int_0^{+\infty} \chi'(-t) t^n Cap_\W(\f<-t) dt
<+\infty, \\
\end{eqnarray*}
This shows that $\f \in \m E_\chi (\W)$.
The other inclusion goes similarly, using the second inequality
in Corollary \ref{est}

Observe that
${\mathcal E}_{\hat{\chi}}(\W) \subset \hat{{\mathcal E}}_{\chi}(\W), \text{
with }  \hat{\chi}(t) = \chi(2t),$ 
as follows by applying 
inequalities of  Corollary \ref{est} 
 with $t=s.$
 
\end{proof}

 Observe that ${\mathcal E}_{\hat{\chi}}(\W) = {\mathcal E}_{\chi}(\W)$ when $\chi(t)=-(-t)^p$.  We thus obtain a characterization
of U.Cegrell's classes ${\mathcal E}^p(\Omega)$
in terms of the speed of decreasing of the capacity of sublevel
sets. This is quite useful since this second definition does not use the
Monge-Amp\`ere measure of the function (nor of its approximants):

\begin{cor}
$$
{\mathcal E}^p(\Omega)=
\left\{ \f \in PSH^-(\W) \, / \,
\int_0^{+\infty} t^{n+p-1} {Cap}_\W(\{\f<-t\}) dt <+\infty \right\}.
$$
\end{cor}

This also provide us with a characterization of the class
${\mathcal F}_a(\Omega)$:

\begin{cor}
$$
{\mathcal F}_a(\Omega)=\bigcup_{\substack{\chi(0) \neq 0,\\ \chi(-\infty)=-\infty}}
{\mathcal E}_{\chi}(\Omega).
$$
\end{cor}

As we shall see in the proof, it is sufficient to consider here
functions $\chi$ that are {\it convex}.

\begin{proof}
The inclusion $\supset $ follows from Proposition \ref{prop}. 
To prove the reverse inclusion, it suffices to show that
if $u\in \mF_a (\W)$ then there exists a function $\chi $ such that 
$u\in \hat{\m E} _\chi(\Omega)$: this is because $\cup \m E_\chi  =\cup  \hat{\m E}_\chi$.
Set
 $$
h(t):= t^n \mrm{Cap}_\W(\{u<-t\}) \ \text{and} \ \tilde{h}(t):=\sup_{s>t} h(s)\ ,\ t>0
$$
 The function $\tilde{h}$ is bounded, decreasing and 
converges to zero at infinity. Consider 
 $\chi (t) := \frac{-1}{\sqrt{ \tilde{h}(-t)}}$ for all $t<0.$ 
Thus $\chi:\R^- \rightarrow \R^-$ is convex increasing, with $\chi(0) \neq 0$
and $\chi(-\infty)=-\infty$. Moreover
\begin{equation*}
\int_0^{+\infty} t^n \chi'(-t) \mrm {Cap}_\W(\{\f<-t\})
dt \le  \frac{1}{2} \int_0 ^{+\infty} \frac{ - \tilde{h}^\prime(s)}{{\tilde{h}^{1/2}(s)}} ds = 
\tilde{h}^{1/2}(0)<+\infty,
\end{equation*}
as follows from  Corollary \ref{est}.
\end{proof}

Let us observe that a negative psh function $u$ belongs to $\mF(\W)$ if and only if 
$ \tilde{h}(0) <+\infty $ (see Corollary 2.5).
\vskip.2cm

We end up this section with the following useful observation.
Let $\chi: \R^- \rightarrow \R^-$ be a non-constant concave increasing function.
Its inverse function $\chi^{-1}:\R^- \rightarrow \R^-$ is convex, hence
for all $\f \in PSH(\W)$, the function
$\chi^{-1} \circ \f$ is plurisubharmonic,
$$
dd^c \chi^{-1} \circ \f =(\chi^{-1})' \circ \f \, dd^c \f+(\chi^{-1})'' d\f \wedge d^c \f
\geq 0.
$$
Now 
$$
Cap_\W(\{\chi^{-1} \circ \f <-t\})=Cap_\W\left(\{\f<\chi(-t) \}\right)
$$
decreases (very) fast if $\chi$ has (very) fast growth at infinity.
Thus $\chi^{-1} \circ \f$ belongs to some class
${\mathcal E}_{\hat{\chi}}(\Omega)$, where $\hat{\chi}$ is completely determined by
$\chi$ and has approximately the same growth order.
This shows in particular that the class
${\mathcal E}_{\chi}(\W)$ characterizes pluripolar sets, whatever the growth of $\chi $:

\begin{thm}
Let $P \subset \Omega$ be a (locally) pluripolar set. Then for any  concave increasing function $\chi:\R^- \rightarrow \R^-$ with $\chi (- \infty) = - \infty$, there exists
$\f \in {\mathcal E}_{{\chi}}(\Omega)$ such that
$$
P \subset \{ \f=-\infty\}.
$$
In particular we can choose $\f \in \mathcal E_{exp} (\Omega),$ where
$$ \mathcal E_{exp} (\Omega) := \left \{ \f \in \mathcal F (\W) ; \int_{\W} e^{- \f} (dd^c \f)^n < + \infty\right\}.$$
\end{thm}

\section{The range of the complex Monge-Amp\`ere operator}

Throughout this section, $\mu$ denotes a fixed positive Borel measure
of finite total mass $\mu(\Omega)<+\infty$  which is dominated by the Monge-Amp\`ere capacity. We want to solve the
following Monge-Amp\`ere equation
$$
(dd^c \f)^n=\mu ,\;
\text{ with } \f \in \m F(\Omega),
$$
and measure how far the (unique) solution $\f$ is from being bounded, by
assuming that $\mu$ is suitable dominated by the Monge-Amp\`ere capacity.

Measures dominated by the Monge-Amp\`ere capacity have been extensively studied by
S.Kolodziej in [K 1,2,3]. The main result of his study, achieved
in [K 2], can be formulated as follows. Fix $\e:\R \rightarrow  [0 ,
\infty [$ a continuous decreasing  function and set 
$F_{\e}(x):=x [\e(-\ln x/n)]^{n}$.
If for all compact subsets $K \subset \Omega$,
$$
\mu(K) \leq F_{\e}(Cap_{\Omega}(K)), \text{ and } \int^{+\infty}
{\e(t)}dt <+\infty,
$$
then $\mu=(dd^c \f)^n$ for some {\it continuous} function
$\f \in PSH(\Omega)$ with $\f_{|\partial \Omega}=0$.

\vskip.1cm

The condition $\int^{+\infty} {\e(t)}dt  <+\infty$ means
that $\e$ decreases fast enough towards zero at infinity. This gives
a quantitative estimate on how fast $\e( -\ln Cap_{\Omega}(K)/n)$,
hence $\mu(K)$, decreases towards zero as $Cap_{\Omega}(K) \rightarrow 0$.

When $\int^{+\infty} \e(t)dt=+\infty$, it is still possible to
show that $\mu=(dd^c \f)^n$ for some function 
$\f \in {\mathcal F}(\Omega)$, but $\f$ will generally be unbounded. 
We now measure how far it is from being so:

\begin{thm} 
Assume for all compact subsets $K  \subset \W$,
\begin{equation} \label{dom}
\mu(K) \leq F_{\e} \Big(\mrm{Cap }_\W(K)\Big ).
\end{equation}

Then there exists a unique function $\f \in \mF (\W)$ 
such that $\mu=(dd^c \f)^n$,
and
\begin{equation*}
\mrm{Cap}_\W(\{\f<-s \}) \leq  \exp (-nH^{-1}(s)), \text{ for all }\ s>0,
\end{equation*}
Here $H^{-1}$ is the reciprocal function of
$H(x) =e\int_{0}^x  \e(t) dt + e \varepsilon (0) + \mu(\W) ^{1/n}$. 

In particular $\f \in \mE _\chi (\Omega)$ with
$-\chi(-t)=\exp (  n H^{-1}(t)/2)$.
\end{thm}

For examples showing that these estimates are essentially sharp, we refer
the reader to section 4 in [BGZ].

\begin{proof} The assumption on $\mu $ implies in particular that it
 vanishes on pluripolar sets.   It follows from \cite{Ce 2} that there exists a unique
$\f \in \mF_a(\W) $ such that  $(dd^c \f )^n = \mu$.
Set
 $$
f(s) : = -\frac{1}{n} \log Cap_\W(\{\f<-s\}),  \ \ \forall s>0
.$$
The function $f$ is increasing  and
$ f(+\infty)=+\infty$, since $Cap _\W$ vanishes on pluripolar sets.

It follows from  Corollary 2.5 and $(5.1)$ that for all $s>0$ and $  t >0$,
$$
t^n Cap_{\W}(\f<-s-t) \leq \mu(\f<-s) \leq F_{\e}\left(Cap_\W(\{\f<-s\})\right).
$$
Therefore
\begin{equation}
\log t - \log \e \circ f (s) +f(s) \leq f(s+t).
\label{estcap}
\end{equation}

We define an increasing sequence $(s_j)_{j \in \N}$ by induction.
Setting
$$ s_{j+1} = s_{j} + e \e \circ f(s_{j}), \text{ for all } j \in \N.$$

\noindent {\it The choice  of $s_0$}.
We choose $s_0\ge 0 $ large enough so that $f(s_0)\ge 0$.
We must insure that $s_0 = s_0 ( \mu ) $ can chosen to be independent of $\f.$ It follows 
from Corollary 2.5 that
\begin{equation*}
\mrm{Cap}_\W(\{\f<-s \}) \le \frac{ \mu(\W)}{s^n}, \ \forall s>0
\end{equation*}
hence  $f(s) \ge \log s  - 1/n \log \mu (\W)$. Therefore $f(s_0)\ge 0 $ if  $s_0= \mu(\W) ^{1/n}$.

\vskip.2cm 
\noindent  {\it The growth of $s_j$. } We can now apply (\ref{estcap}) and get
 $f(s_j) \ge j + f(s_0) \ge j.$
Thus  $\lim_j f(s_j)=+\infty$.
There are two cases to be considered. 

If $s_{\infty}=\lim s_j \in \R^+$, then
$f(s) \equiv +\infty$ for $s > s_\infty$, i.e. 
$ Cap_\W(\f<-s)=0,  \ \ \forall s  > s_\infty$. 
Therefore $\f$ is bounded from below by $-s_\infty$,
in particular $\f \in \mE_\chi(\W)$ for all $\chi .$

Assume now ( second case) that $s_j\to +\infty.$ For each   
 $s>0, $  there exists $ N = N_s  \in \N $ such that
 $s_{N} \le s < s_{N +1}.$
We can estimate $s \mapsto N_s$,
\begin{eqnarray*}
  s \le s_{N +1}&=&  \sum_0^{N} (s_{j+1} -s_j) + s_0  = \sum_0^{N}
   e \,  \e \circ f(s_{j}) + s_0 \\
&  \le &  e \sum_{0}^{N}
 \e (j) + s_0 
\le e \int_{0}^{N}
  \e (t)dt + \tilde{s}_0  =:   H(N) ,
\end{eqnarray*}
where $  \tilde{s}_0  = s_0 + e.\e(0).$ 
Therefore $H^{-1}(s) \le N \le f(s_N) \le f(s),$ hence
\begin{equation*}
Cap_\W(\f<-s) \leq  \exp (-nH^{-1}(s)).
\end{equation*}

Set now $ g(t) = -\chi(-t) =\exp (  n  H^{-1}(t)/2)$. Then
\begin{multline*}
\int_0^{+\infty} t^n g'(t) Cap_\W(\f<-t) dt\\  \le
\frac{n}{2}\int_0^{+\infty} {  t^n}\frac{1}{\e(  H^{-1}(t)) + s_0 } \exp (-n H^{-1}(t)/2) dt\\
\le C \int_{0}^{+\infty}  { (t+1)^n} \exp (n(\alpha-1) t ) dt <+\infty.
\end{multline*}
This shows that  $\f \in \m E_\chi (\W)$ where $\chi(t)=-\exp (  n  H^{-1}(-t)/2)$.
\end{proof}
Observe that the proof above gives easily an a priori uniform bound of the solution of $(dd^c \f)^n = \mu$, when $\mu$ is a finite Borel mesure on $\W$ satisfying $(5.1)$ with  $\int_0^{+ \infty} \varepsilon (t) d t < + \infty$ (see also [K2]). Indeed it follows from the above estimates that $\f \geq - s_{\infty},$ where 
$$s_{\infty}  \leq e \int_0^{+ \infty} \varepsilon (t) d t + e \varepsilon (0) + \mu (\Omega)^{1 \slash n}.$$

We now generalize U.Cegrell's main result [Ce 1].

\begin{thm} Let $\chi : \ \R^- \to \R^- $ be
an increasing function such that $\chi(-\infty) = -\infty $.
Suppose there exists a locally bounded function $ F : \ \R^+ \to \R^+ $ 
such that $\limsup_{t \to +\infty} F(t)/t<1,$ and 
\begin{equation} \label{quant}
\int _\W (-\vk) \circ u \, d\mu \le F( E_\vk (u )), \ \ \forall \ u \in
\m T (\W ),
\end{equation}
where $E_{\chi}(u):=\int_{\Omega} (-\chi) \circ u (dd^c u)^n$
denotes the $\chi$-energy of $u$.

Then there exists a  function  $ \f \in \mE_\vk (\W )$ 
such that $\mu = (dd^c \f)^n.$
\end{thm}

\begin{proof}
 The assumption on $\mu $ implies in particular
 that it vanishes on pluripolar sets. It follows from \cite{Ce 2}
 that there exists a function $u \in \m T(\W)$ and
 $f \in L_{loc}^1 \big((dd^c u )^n\big)$ such that 
$\mu = f (dd^c u )^n. $ 

Consider $\mu_j:=\min (f, j ) (dd^c u )^n$.
This is a finite measure which is bounded from above by 
the Monge-Amp\`ere measure of a bounded function. It follows therefore from
\cite{K 1} that there exist $\f_j \in \m T (\W)$ such that 
$$
(dd^c \f_j )^n = \min (f, j ) (dd^c u )^n.
$$
 The comparison principle shows that $\f_j$ is a decreasing sequence.
Set  $\f =\lim_{j\to \infty } \f_j$.
It follows from (\ref{quant}) that
$E_\vk (\f_j) ( F(E_\vk (\f_j)))^{-1} \le 1 $, hence
$\sup_{j\ge 1} E_\vk (\f_j) <\infty . $
This yields $\f \in \mE_\vk(\Omega) $.

We conclude now by continuity of the Monge-Amp\`ere operator
along decreasing sequences that
 $(dd^c \f )^n = \mu.$ 
\end{proof}

When $\chi(t)=-(-t)^p$ (class ${\mathcal F}^p(\Omega)$), $p \geq 1$,
the above result was established by U.Cegrell in [Ce 1].
Condition (5.3) is also necessary in this case,
and the function $F$ can be made quite explicit: there exists
$\f \in {\mathcal F}^p(\Omega)$ such that $\mu=(dd^c \f)^n$
if and only if $\mu$ satisfies (5.3) with $F(t)=C t^{p/(p+n)}$,
for some constant $C>0$. 

Actually the measure  $\mu$ satisfies (5.3) for $\chi(t)=-(-t)^p$, and  $F(t)=C\cdot t^{p/(p+n)}, $ $p>0 $  if and only if  $\m F^p(\W) \subset L^p(\mu) $ (see [GZ]).
 
We finally remark that this condition  can be interpreted in terms of domination by capacity.

\begin{pro}
If $\m F^p(\W) \subset L^p(\mu) ,$ 
then there exists $C>0$ such that
$$
\mu(K) \leq C \cdot  Cap_{\Omega}(K)^{\frac{p}{p+n}},
\; \text{ for all } K \subset \Omega.
$$
Conversely if $\mu(\cdot) \lesssim Cap_{\Omega}^{\a}(\cdot)$ for some $\a>p/(p+n)$, then
$\m F^p(\W) \subset L^p(\mu) .$ 
\end{pro}

\begin{proof}
 The estimate   (\ref{quant}) applied  to $u = u_K ^* $, the relative extremal function of the compact $K$, yields  
\begin{eqnarray*}
\mu (K)& = &\int_\W \mrm 1 _K \cdot  d\mu  \le  \int_\W (-u^* _K)^p d\mu \\
 &\le &  C\cdot \left( \int_\W (- u^*_K)^{ p} (dd^c u_K^*)^n \right ) ^\frac{p}{p+n}\\
&=& C\cdot \left [ Cap_\W (K) \right ]^\frac{p}{n+p}.
\end{eqnarray*}

Conversely, assume that   $\mu(K) \le C.  Cap_{\Omega}^{\a}(K)$ for all compact $K\subset \W$, where $\a > p/(n+p)$ then (\ref{quant}) is satisfied.
Indeed, if  $u \in \m F^p(\W) ,$ then 
\begin{multline*}
\int_\W (-u)^p d\mu = p\int _1 ^\infty t^{p-1} \mu ( u<-t) dt + O (1)\\
 \le  C\cdot p \int _1 ^\infty t^{p-1} \big(Cap_\W ( u<-t) \big )^\a  dt + O (1)\\
\le   C\cdot \Big( \int _1 ^\infty t^{n+p-1} Cap_\W ( u<-t)dt \Big )^\a  \cdot 
\Big (\int_1^\infty t^{[p-1-\a (n+p-1)]/\beta}  dt \Big )^\beta + O(1) ,
\end{multline*} 
where $\a + \beta =1. $ 
The first integral converges by Corollary 4.3, the latter one is finite since $p -1 - \a (n+p-1) > \a - 1 = - \beta.$
\end{proof}

\vskip .2cm

Slimane Benelkourchi, Vincent Guedj and Ahmed Zeriahi

Institut de Math\'ematiques de Toulouse,

Laboratoire Emile Picard, 

Universit\'e Paul Sabatier

118 route de Narbonne

31062 TOULOUSE Cedex 09 (FRANCE)

benel@math.ups-tlse.fr

guedj@math.ups-tlse.fr

zeriahi@math.ups-tlse.fr

\end{document}